\documentclass[12pt,a4paper]{article}
\usepackage{amscd}
\usepackage{amsmath}
\usepackage{amsfonts}
\usepackage{amssymb}

\usepackage[dvips]{color}

\input{epsf}             
\def\epsfcenter#1{{\vcenter{\hbox{\epsfbox{#1}}}}} 
\def\epsfpad#1#2#3{{\vcenter{\vskip#2\hbox{\epsfbox{#1}}\vskip#3}}} 

\def \S {\Sigma}
\def \n {\widetilde}

\newcommand{\CH}{{\rm CH}} 

\newcommand{\sll}{{\rm sl}}   

\newcommand{\U}{{\rm U}}

\DeclareMathOperator{\interior}{int}

\newcommand{\C}{{\mathbb C}} 
\newcommand{\Z}{{\mathbb Z}} 
 
\newtheorem{theorem}{Theorem}
\newtheorem{lemma}{Lemma}

\newtheorem{definition}{Definition}

\begin{document}

\title{ Observables in the Turaev-Viro and Crane-Yetter models}  
\author{John W.\ Barrett\\
\small{\it School of Mathematical Sciences, University Park, Nottingham NG7 2RD, UK} \and
J. Manuel Garc\'\i a-Islas\\
\small{\it Centro de Investigacion en Matem\'aticas,
A.P. 402, 36000, Guanajuato, Gto, Mexico}
\and
Jo\~ao Faria Martins \footnote{Also at Departamento de Matem\'atica, Universidade Lus\'{o}fona de
Humanidades e Tecnologia, Av do Campo Grande, 376, 1749-024, Lisboa,
Portugal. }
\\
\small{\it Departamento de Matem\'atica,
Instituto Superior T\'ecnico,}\\
\small{\it Av. Rovisco Pais,
1049-001 Lisboa,
Portugal}
  }
 \date{}
 
 \maketitle 

\begin{abstract}   We define an invariant of graphs embedded in a three-manifold and a partition function for 2-complexes embedded in a triangulated four-manifold by specifying the values of variables in the Turaev-Viro and Crane-Yetter state sum models. In the case of the three-dimensional invariant, we prove a duality formula relating its Fourier transform to another invariant defined via the coloured Jones polynomial. In the case of the four-dimensional partition function, we give a formula for it in terms of a regular neighbourhood of the 2-complex and the signature of its complement. Some examples are computed which show that the partition function determines an invariant which can detect non locally-flat surfaces in a four-manifold.
\end{abstract}

\section{Introduction} 
In this paper we define an invariant of graphs embedded in a three-manifold and a partition function for 2-complexes embedded in a triangulated four-manifold.   The invariance under homeomorphisms\footnote{In the main results, homeomorphism means PL homeomorphism. However an alternative proof of invariance is also sketched for smooth manifolds and maps.} of the manifold is proved in the case of graphs embedded in a three-manifold. It is shown how the invariant is determined by the Fourier transform of a previously defined invariant, the relativistic spin network invariant, which is defined using the technology derived from the Jones polynomial of links.  This provides a proof of the general result announced in \cite{BAR} which was only proved there in simple cases. In the case of a 2-complex embedded in a four-manifold, we give a formula for the partition function which involves only a regular neighbourhood of the 2-complex. These formulae enable a simple evaluation in the case where the 2-complex is a surface.
We conjecture that the partition function can be normalised to give a
homeomorphism invariant of the 2-complex in the 4-manifold. This is proved in a particular special case.

The invariants described here are important in the study of quantum gravity. The 3-manifold invariant describes the partition function for three-dimensional quantum gravity in the presence of the observation of distance or momentum observables \cite{BAR, BFY}. In four dimensional quantum gravity, it has been suggested that there should be a perturbative expansion of quantum gravity about a topologically-invariant vacuum \cite{SS,FS}. Each term in the expansion would be determined by the expectation value of some observables in a topologically-invariant partition function. This idea motivates the definition of the four-dimensional partition function presented here.

Let $M$ be a triangulated closed 3-manifold and $r\ge 3$ an integer. Then the
Turaev-Viro invariant of $M$ is a homeomorphism invariant of $M$ defined using a
state-sum model on the triangulation of $M$ \cite{TV}. Each edge in the triangulation
is labelled with a half-integer variable $j$, the spin, from the set
${\cal L}=\{0,1/2,1,\ldots, (r-2)/2 \}$. The statistical weights for each simplex are
constructed so that after summing over the variable $j$ for each edge, the
resulting partition function $Z(M)\in\C$ is independent of the triangulation of $M$. The full details of this definition of the Turaev-Viro invariant are not given here; however an alternative description of the invariant in terms of the chain-mail link is given below.

Now let $\Gamma$ be a graph determined by a subset of the set of edges of the triangulation and label each edge in $\Gamma$ with a spin from $\cal L$. This labelled graph in $M$ is denoted $\Gamma(j_1,j_2,\ldots,j_n)$.
Define the partition function $Z(M,\Gamma(j_1,j_2,\ldots,j_n))$ in the same way as the Turaev-Viro partition function but not summing over the variables $j_i$. We show that this definition is actually a homeomorphism invariant of $(M,\Gamma)$, that is, does not depend on the triangulation of $M$ away from $\Gamma$.
In fact the definition can be extended to give a partition function for any labelled graph in $M$, not just those that can be realised as edges of some triangulation. 

As a special case one can recover the refined invariant of Turaev and Viro
 \cite{TV,ROEX,YEX}, which depends on a choice of an element of $H^1(M,\Z_2)$, by taking $\Gamma$ to be the 1-skeleton of a triangulation and summing over spins on $\Gamma$ whose parity agrees with a representative cocycle for the element of $H^1(M,\Z_2)$. 
 
The following is the main result about this partition function.
It relates the partition function to the relativistic spin network invariant $Z_R(M,\Gamma(k_1,k_2,\ldots,k_n))$ of the labelled graph $\Gamma$ in a connected closed orientable 3-manifold $M$. This relativistic invariant was previously defined by Yokota and Yetter for the case $M=S^3$ \cite{YO,Y,BCLA}, and was used in 4-dimensional quantum gravity \cite{BC}. A generalisation to any $M$ is defined below by a suitable squaring of the Witten-Reshetikhin-Turaev invariant of a graph embedded in a 3-manifold \cite{RT}. 

The result requires a matrix which is the kernel of the Fourier transform.
Let $H_{jk}$ be the matrix determined by the coloured Jones polynomial for the Hopf link, and $\dim_q j$ the coloured Jones polynomial for the unknot. Explicitly,
$$H_{jk} =  \frac{\sin\frac\pi r(2j+1)(2k+1)}{\sin\frac\pi r}(-1)^{2j+2k}$$ and
$$\dim_qj=H_{0j}=\frac{\sin\frac\pi r(2j+1)}{\sin\frac\pi r}(-1)^{2j}.$$

\begin{theorem}[Fourier transform] \label{theoremft}For any triangulation of $M$,
\begin{multline}
\sum _{{j_1j_2 \dots}  j_n} 
{Z(M,\Gamma(j_1,j_2,\ldots,j_n))} 
\frac{H_{j_1 k_1}}{\dim_qj_1} 
\frac{H_{j_2k_2}}{\dim_qj_2}\dots \frac{H_{j_n k_n}}{\dim_qj_n}\\
=   Z_R(M,\Gamma(k_1,k_2,\ldots,k_n)).
\end{multline}
 \end{theorem}  
 
In section \ref{3dInvariance} it is shown how to extend these definitions and results to an arbitrary graph embedded in $M$.
 
Now we describe the formula for the corresponding 4-manifold construction. Let $W$ be a closed triangulated 
 oriented 4-manifold. The Crane-Yetter state-sum invariant of $W$ is 
 given by assigning a spin variable to each 2-simplex of the triangulation, and summing over all values of these variables with weights determined by the sums of products of 15-j symbols \cite{CKY}. Now suppose $\Gamma$ is a 2-complex determined by a subset of the 2-simplexes of the triangulation, and label each 2-simplex with a spin from $\cal L$. The labelled 2-complex is denoted $\Gamma(j_1,j_2,\ldots,j_n)$, with $n$ the number of 2-simplexes, and again the partition function $Z(W,\Gamma(j_1,j_2,\ldots,j_n))$ is determined in the same way as the Crane-Yetter partition function, but not summing over the variables $j_i$.
 
The main result is a formula for this partition function in terms of the Witten-Reshetikhin-Turaev invariant $Z_{WRT}$ of a labelled link in a 3-manifold. 
A regular neighbourhood of $\Gamma$, denoted $\widetilde\Gamma$, is a manifold with boundary. 
Its interior is denoted $\interior(\widetilde\Gamma)$, and the complement $W\setminus\interior(\widetilde\Gamma)$, a compact manifold with boundary.
In this boundary there is a framed link $\gamma$ determined by a circumference for each 2-simplex of $\Gamma$  (or attaching curve for the corresponding 2-handle). Each component of the link $\gamma_i$ is labelled with the corresponding spin $j_i$. 
Let $\kappa=\exp(i\pi(r-2)(3-2r)/4r)$, $N=\sum_j(\dim_q j)^2$, $\sigma(W)$ the signature (index) of $W$ and $\chi(\Gamma)$ the Euler characteristic of $\Gamma$.

The result for the partition function is

\begin{theorem}[4d formula]\label{theoremfd} For any triangulation of $W$,
\begin{multline} 
Z\bigl(W,\Gamma(j_1,j_2,\ldots,j_n)\bigr)
 =  \\ 
 Z_{WRT}\bigl(\partial(W\setminus\interior(\widetilde\Gamma)),\gamma(j_1,j_2,\ldots,j_n)\bigr) \kappa^{\sigma(W\setminus\interior(\widetilde\Gamma))} N^{\frac{\chi(\Gamma)}2-n}\dim_qj_1\ldots\dim_qj_n.
\end{multline}
 \end{theorem}  
Note that the 3-manifold  $\partial(W\setminus\interior(\widetilde\Gamma))$ is the same as $\partial(\widetilde\Gamma)$ but with the opposite orientation. The importance of expressing the observable in this way is that $Z_{WRT}$ is known to be an invariant of the link $\gamma$ in the 3-manifold. This leads to the easy computation of simple examples. 

Section \ref{sectionexamples} computes some examples which demonstrate the dependence of this definition  of  the partition function on  the triangulations of the 2-complex and of the 4-manifold. In the case in which all the faces are assigned the trivial representation,  then  the result is a homeomorphism invariant of the 2-complex embedded in the 4-manifold.  This invariant is non-trivial in cases of non locally-flat submanifolds.  

 The precise definitions and the proof of the theorems appear below.

\section{Chain-mail}\label{secchm}
The main technique is the description of the Turaev-Viro invariant called chain-mail\cite{RO}.\footnote{The name chain-mail comes from the technique of making body armour with interlaced metal rings.}  This starts from a handle decomposition of the manifold, rather than a triangulation, and associates a framed link in $S^3$ (the chain-mail link) to each handle decomposition. The Turaev-Viro partition function can then be evaluated from the coloured Jones polynomial of this link. This is carried out for $Z(M)$ in \cite{RO} and the description given here is an extension of this technique to the case $Z(M,\Gamma)$ considered here. Each triangulation of the manifold determines a canonical handle decomposition, so there is no loss of generality. Handle decompositions are considerably more flexible than triangulations and the proofs appear to be easier to construct using handles. There is the added bonus that the invariant is defined for more general graphs than those that occur as subsets of edges in a triangulation (see section \ref{3dInvariance}).

\subsection{The chain-mail link}
 The 3-manifolds in this section are required to be orientable.
Let $M$ be a closed compact 3-manifold with a decomposition into handles.
Let $H_+$ be the union of the 0-handles and 1-handles, and $H_-$ the union of the 2- and 3-handles. This gives a Heegaard splitting of $M$, a decomposition into two  oriented handlebodies. A generalised Heegaard diagram for the handle decomposition of $M$ consists of the handlebody $H_-$ together with the set $\{\gamma_+\}$ of attaching curves in $\partial H_-$ for the meridianal disks of the 1-handles, and the set $\{\gamma_-\}$ of attaching curves in $\partial H_-$  for the 2-handles. Each set consists of a non-intersecting set of simple closed curves, though of course the curves in $\{\gamma_+\}$ can intersect the curves in $\{\gamma_-\}$.
  
Given the Heegaard splitting and its diagram $D=(H_-,\{\gamma_+\},\{\gamma_-\})$, the chain-mail link in $S^3$ is defined as follows. Pick an  orientation preserving embedding $\Phi\colon H_-\to S^3$. The components of the chain-mail link consist of the $\Phi(\gamma_-)$ and the $\Phi(\overline\gamma_+)$, where the curves $\overline\gamma_+$ are displaced into the interior of $H_-$ in a level surface of a collar neighbourhood of the boundary. Moreover, the link inherits a framing by considering each curve $\gamma_\pm$ to be a thin strip parallel to the boundary of $H_-$. This defines the framed link $\CH(D,\Phi)$. 
  
\subsection{The chain-mail invariant}

Let $L$ be a framed link in $S^3$.  Suppose that each component of the link is assigned a spin. The labelled link is denoted $L[j_1,j_2,\ldots,j_K]$, where $K$ is the total number of components and $j_i$ is the spin of the 
$i$-th component. The coloured Jones polynomial gives an evaluation of the link  $<L[j_1,\ldots,j_K]>\;\in\C$ for each parameter $r$. 
The evaluation is defined by projecting the link to a two-dimensional diagram and using the Kauffman bracket evaluation of the link diagram in which the $i$-th component of the link is cabled $2j_i$ times and the Jones-Wenzl idempotent corresponding to an irreducible representation of quantum $\sll_2$ is inserted into the cable. This evaluation is described in \cite{KL}; we use the convention that the empty link has evaluation 1, and Kauffman's parameter $A$ is given by $A=e^{\frac{i\pi}{2r}}$. The evaluation does not depend on the choice of the projection.

Each 1-handle of the manifold $M$ corresponds to exactly one circle $\gamma_+$ of $D$ and each 2-handle to exactly one $\gamma_-$ of $D$. Let $h_1$ be the number of 1-handles and $h_2$ the number of 2-handles. Then the chain-mail link $\CH(D,\Phi)$ has $K=h_1+h_2$ components. 
Now let $\widetilde\Gamma$ be a submanifold consisting of a subset of the 0- and 1-handles (so as to include all the 0-handles that the 1-handles attach to), with a spin label $j_i$ for the $i$-th 1-handle in $\widetilde\Gamma$.
Assume that the first $n$ components of the link correspond to the $n$ 1-handles of $\widetilde\Gamma$. Let $g$ be the genus of the Heegaard surface, and $N$ the constant $N=\sum_j(\dim_q j)^2$.
Then the chain-mail version of the partition function is given by

\begin{definition} \label{chainmail}
\begin{multline*}
 Z_{\CH}(M,\widetilde\Gamma(j_1,j_2,\ldots,j_n))= \\
 N^{g-K-1} \sum _{{j_{n+1}j_{n+2} \dots}  j_K}
  \left<\CH(D,\Phi)[j_1,\ldots,j_K]\right>\;\prod_{i=1}^K \dim_q j_i 
\end{multline*}
\end{definition}  
  
Two results about this definition follow immediately.

\begin{lemma} $Z_{\CH}$ is independent of $\Phi$.\end{lemma}

A triangulation of $M$ induces a handle decomposition of $M$ by thickening the simplexes in $M$ (so that the $k$-skeleton is thickened to a regular neighbourhood of itself).

\begin{lemma} \label{lemmaind} $Z(M,\Gamma(j_1,j_2,\ldots,j_n))=Z_{\CH}(M,\widetilde\Gamma(j_1,j_2,\ldots,j_n))$.\end{lemma}

 The proof of these two lemmas is identical to that given by Roberts for the case where $\Gamma$ is empty. Each circle $\gamma_-$ `kills' the spin passing through it. This means that the embedding of the 1-handles of $H_-$ (2-handles of $M$) can be rearranged. In the second case, this leaves a product of tetrahedral graph evaluations (6$j$-symbols) which defines the Turaev-Viro partition function for a triangulation.
 
\subsection{Invariance}

The following result shows that the partition function with non-trivial $\Gamma$ does not depend on the triangulation of $M$ away from $\Gamma$ itself. 

\begin{theorem}[Invariance] \label{theoreminvt}Let $M'$ be a second triangulated 3-manifold with labelled subset of edges $\Gamma'$ and $\phi \colon M\to M'$ a homeomorphism such that $\phi(\Gamma)=\Gamma'$ as labelled graphs. Then
$$Z(M,\Gamma(j_1,j_2,\ldots,j_n))=Z(M',\Gamma'(j_1,j_2,\ldots,j_n)).$$
\end{theorem}

This theorem follows immediately from theorem \ref{theoremft}, proved below, and the invariance of the relativistic spin network invariant, since the Fourier transform is invertible (see section \ref{3dInvariance}). However we sketch a direct proof in the differentiable category which does not involve the use of theorem \ref{theoremft}. 

\noindent{\it Proof.}
To establish this, we begin with the corresponding property for $Z_{\CH}$.
Consider two different handle decompositions of $M$ which agree on 
$\widetilde\Gamma$. A handle decomposition corresponds to a Morse function on $M$. Any two Morse functions can be deformed one into the other, inducing a set of moves on the corresponding generalised Heegaard diagram \cite{KU}. In our case, the Morse functions can be arranged so that $\partial\widetilde\Gamma$ is a level surface, and the values on $\widetilde\Gamma$ are less than the values on the remainder of $M$ and do not vary when the Morse function is deformed. The consequence of this is that in the moves on the generalised Heegaard diagram listed in \cite{KU}, there are no instances of 1-handle curves sliding over the attaching curves for the 1-handles
of $\widetilde\Gamma$. Similarly, the attaching curves for the 1-handles
of $\widetilde\Gamma$ never vanish. It is easily seen that the partition function $Z$ is invariant under these moves.

The triangulated manifolds $M$ and $M'$ both have canonical handle decompositions; although the map $\phi$ does not necessarily take the 1-handles and their attaching curves corresponding to $\Gamma$ to those corresponding to $\Gamma'$, it can be adjusted by an isotopy so that it does. This follows from the fact that these 1-handles and the 0-handles they are attached to form tubular neighbourhoods 
of $\Gamma$ and  $\Gamma'$, and by using standard results on tubular neighbourhoods \cite{KOS}. This completes the proof of theorem \ref{theoreminvt}.

\section{$S^3$ and the relativistic invariant}

In the case of $M=S^3$ there is another known invariant of embedded graphs with edges labelled by spins, the relativistic spin network invariant. In this section it is shown by explicit calculation that the two invariants are related by a Fourier transform of the spin labels. The calculation shows a relation to the shadow world presentation of the coloured Jones polynomial of Kirillov and Reshetikhin \cite{KR,KL}.

We start by giving a definition of the relativistic spin network invariant for a graph embedded in $S^3$ in terms of the Kauffman bracket evaluation of a diagram of the graph in $S^2$. This diagram is obtained by assuming the graph does not meet the north and south poles of $S^3$ (moving it if necessary). These two points can be removed and the resulting space projected to the equator, $S^2$. 

Then we give a process for turning a diagram of the embedded graph into a chain-mail link for $S^3$. This is based on a handle decomposition of $S^3$ which includes the original graph in the 0- and 1-handles. So it is possible to give a description of the partition function for the embedded graph in terms of the evaluation of this chain-mail link. Finally, we show that this chain-mail evaluation is the Fourier transform of the relativistic spin network evaluation for the original graph.

\subsection{Relativistic spin networks} \label{rsn}
The definition of the relativistic spin network invariant is as follows \cite{BCLA}.  It is convenient to generalise the usual definition of graph (which has edges which are closed intervals with two ends which lie on vertices) to allow, in addition, components which are circles with no vertices (these will not be called edges as they are topologically distinct). A graph is required to have a finite number of edges and circles; thus it forms a compact 1-dimensional polyhedron.

Let $\Gamma(i_1,i_2,\ldots,i_n)$
be a  graph embedded in $S^3$, with its edges and circles labelled with spins $i_1,i_2,\ldots,i_n$ (in a fixed order). 
First, the invariant is defined in the case of graphs in which each vertex is trivalent, then this will be generalised to arbitrary vertices.
For each vertex of a trivalent graph there are three spin labels $(i,j,k)$ on the three edges meeting the vertex. 
Put 
$$\Theta(i,j,k)=\left<\epsfcenter{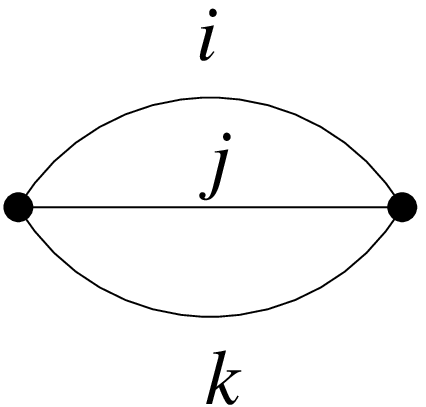}\right>$$
Then the relativistic invariant $\left\langle \Gamma(i_1,i_2,\ldots,i_n) \right\rangle_R$ is defined in terms of
the Kauffman bracket invariant of the diagram given
by projecting the graph in $S^3$ to $S^2$ by
\begin{equation}\label{reldef}
\left\langle \Gamma \right\rangle_R = 
\frac{\left| \left\langle \Gamma \right\rangle \right| ^2} {\prod_{\text{vertices}}\Theta}.
\end{equation}
We note that this definition includes the case of knots and links, for which there are no vertices.

The definition is extended to arbitrary graphs by the relations
\begin{equation*}
\left\langle \epsfcenter{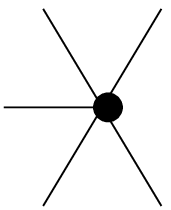}\,\right\rangle _R=\sum _j \left\langle \epsfcenter{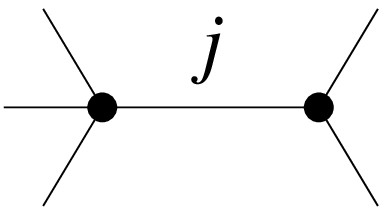} \,\right\rangle _R \
\dim _q j  
\end{equation*}
which defines an $n$-valent vertex recursively, for $n>3$, 
\begin{equation*}
   \left\langle  \, \epsfcenter{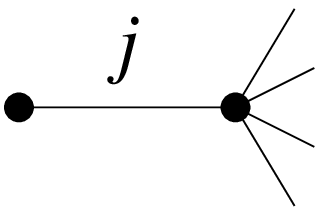}\right\rangle _R=\delta _{j0}
\left\langle \, \epsfcenter{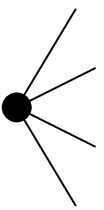} \right\rangle _R 
\end{equation*}
for 1-valent vertices, and
\begin{equation}\label{twovalent}
\left\langle \epsfpad{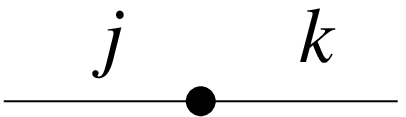}{0pt}{10pt}  \right\rangle _R=\frac{1}{\dim _qj} \delta _{jk}
\left\langle  \epsfpad{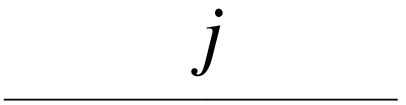}{0pt}{10pt} \right\rangle _R
\end{equation}
for 2-valent vertices. In these equations, only the part of the graph is shown, the remainder being the same on both sides.
By using these relations, an arbitrary graph is reduced to one in which each component is either a circle, or consists of trivalent vertices with edges which are closed intervals. 
It is worth remarking that one can also view this invariant as one determined by a particular subcategory of representations of  $\U_q\sll_2\times\U_{q^{-1}}\sll_2$.

\subsection{Graph diagrams to chain-mail}
The process for turning an embedded graph diagram in $S^2$ into chain-mail is to construct a handle decomposition of $S^3$ which contains $\Gamma$ in its 0- and 1-handles.  Then the handle decomposition is turned into chain-mail as in section \ref{secchm}.  For convenience in presenting the proof, we shall assume that each (two-dimensional) region of the diagram is a topological disk. This implies that the diagram is connected and contains at least one closed circuit. The proof can be easily extended to other cases by adjusting the normalising factors.  Also, we assume that if the diagram is the unknot, then it contains at least one vertex, so that the circle is cut into edges.

The procedure starts by adding extra vertices and edges to the graph $\Gamma$ to make a new graph $\Gamma'$ for which it is possible to fill in 2- and 3-handles to make $S^3$ in an obvious way.  At each crossing in the diagram, add an extra vertex to the graph on both of the edges that cross (so that both new vertices project onto the crossing point), and add a new edge connecting them, vertically with respect to the projection.  
The new graph $\Gamma'$ is labelled in the following way. Where an edge has been subdivided into two edges by the introduction of one of the new vertices, one of the two new edges inherits the same label, while the other new edge is labelled with 0. (It does not matter which is which). The labels on the new connecting edges are summed over, weighted with $\dim_q j$. Let $C$ be the number of crossings in the diagram. With this labelling,
$$ Z(S^3,\Gamma)=N^{2C}
\sum_{
  \genfrac{}{}{0pt}{} {\text{connecting}}{\text{edge labels j}}
}
Z(S^3,\Gamma')\prod\dim_q j$$
This resulting graph $\Gamma'$ can now be thickened to provide the 0- and 1-handles for a handle decomposition of $S^3$. The 2-handles are now just a thickening of the planar regions in the diagram, which are disks by the assumption above. The 0- 1- and 2-handles together give a space homeomorphic to $S^2\times [0,1]$, and so adding two 3-handles, one above and one below the diagram, gives $S^3$. This handle decomposition contains the original graph as a subset of the 0- and 1-handles.

The chain mail diagram  for this handle decomposition can be obtained very simply from the original diagram for $\Gamma$. 
Make the following replacements
$$
\epsfcenter{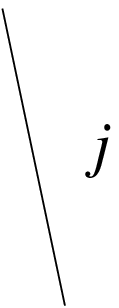}\qquad\to\qquad \epsfcenter{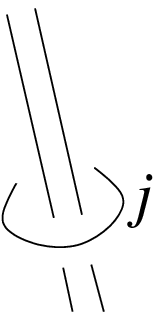}
$$
$$
\epsfcenter{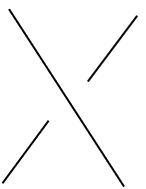}\qquad\to\qquad \epsfcenter{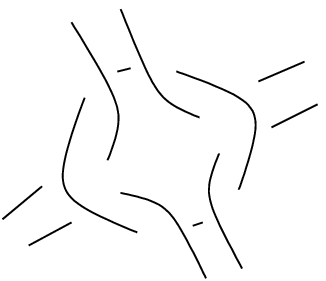}
$$
$$
\epsfcenter{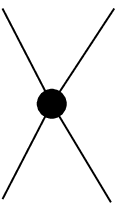}\qquad\to\qquad \epsfcenter{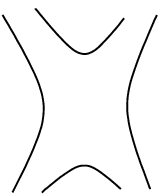}
$$
with a similar replacement to the last one for an $n$-valent vertex, for any $n$.
In the first of these replacements it has to be understood that there is just one circle marked with $j$ for each complete edge of $\Gamma$. If the edge crosses other edges then the circle is placed on one of the resulting segments between two crossings, and it does not matter which segment is chosen.

\noindent{\it Example.} The chain-mail diagram for the labelled graph
$$\epsfcenter{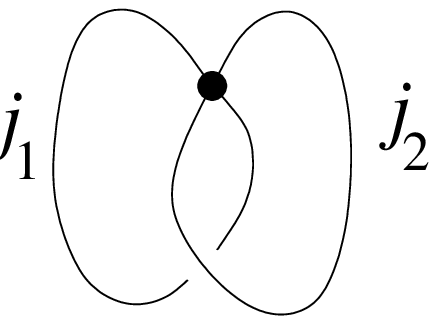} $$
is
$$\epsfcenter{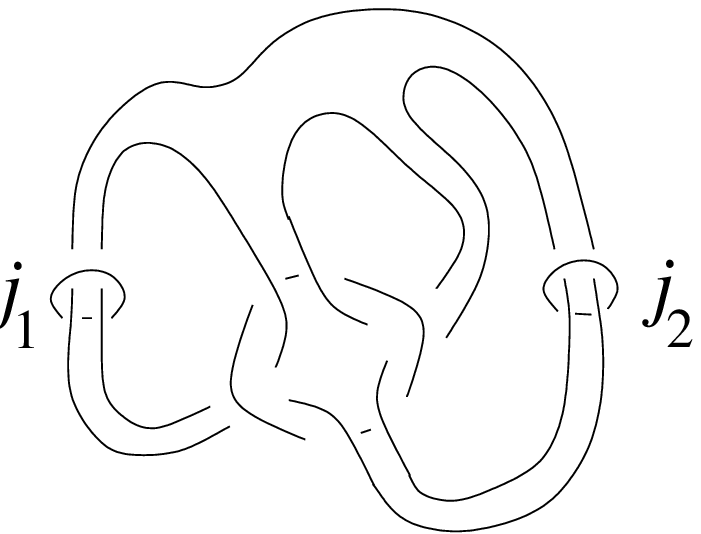}$$
  
In this diagram, the circles corresponding to edges of $\Gamma'$ labelled by $0$ are not shown, as they do not affect the evaluation.
  The complete circle shown in the first two replacements moves are the $\gamma_+$ curves, with the curves resulting from the first replacement move labelled with a spin. The remaining circles formed from the fragments shown in each replacement form the $\gamma_-$ circles.

\subsection{Shadow world evaluation}
The chain-mail diagram is evaluated using definition \ref{chainmail} to give the partition function
$Z(S^3,\Gamma')$. The normalisation factor $N^{g-K-1}$ in definition \ref{chainmail} reduces to $N^{-2-E-3C}$, where $E$ is the number of edges in $\Gamma$ and $C$ the number of crossings in the diagram. This follows from the fact that $-(g-K-1)$ is the number of 1-handles plus the number of 3-handles, and there are $E+3C$ 1-handles and two 3-handles.

Using this explicit presentation of a formula for the partition function, it is easy to give a proof of theorem \ref{theoremft} for the case of $M=S^3$, namely
\begin{theorem}\label{theoremftcase}
\begin{multline}
\sum _{{j_1j_2 \dots}  j_n} 
\frac 
{Z(S^3,\Gamma(j_1,j_2,\ldots,j_n))} 
{Z(S^3)}
\frac{H_{j_1k_1}}{\dim_qj_1} 
\frac{H_{j_2k_2}}{\dim_qj_2}\dots \frac{H_{j_nk_n}}{\dim_qj_n}\\
=   \left<\Gamma(k_1,k_2,\ldots,k_n)\right>_R,
\end{multline}
\end{theorem}
as stated in \cite{BAR}. Examples of such calculations are given in \cite{GI}. In the rest of this section, we outline the proof.
 
Starting with the evaluation of the chain-mail link using definition \ref{chainmail}, the expression can be transformed using the following identities:
$$
\sum_j H_{jk}\left< \epsfcenter{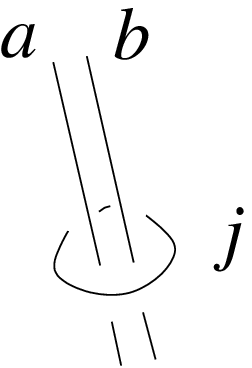}\right> \qquad= \qquad\frac{N}{ \theta_{abk}}\left< \epsfcenter{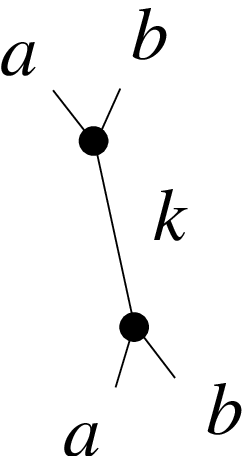}\right>
$$
$$
\left<\epsfcenter{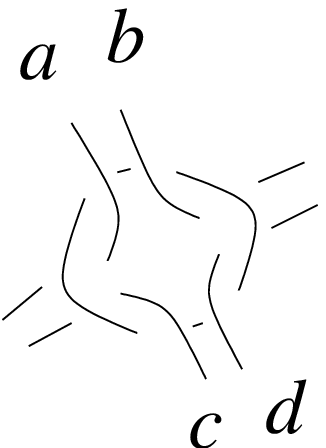}\right> \qquad=\qquad \sum_{ef} \frac{N\dim_qe\,\dim_qf}{\theta_{abe}\theta_{cde}\theta_{acf}\theta_{bdf}}
\left<\epsfcenter{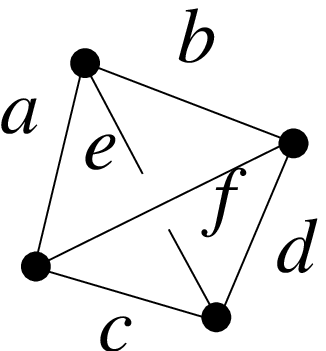}\right> \left<\epsfcenter{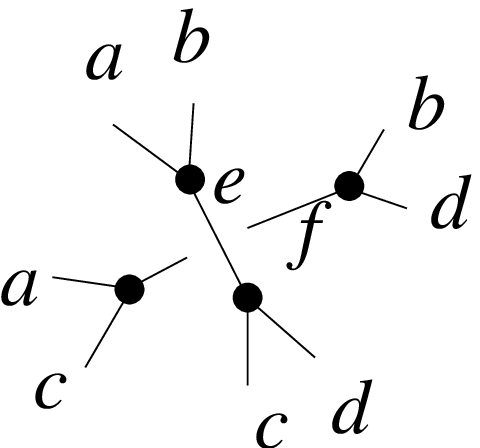}\right>
$$
$$
\left<\epsfcenter{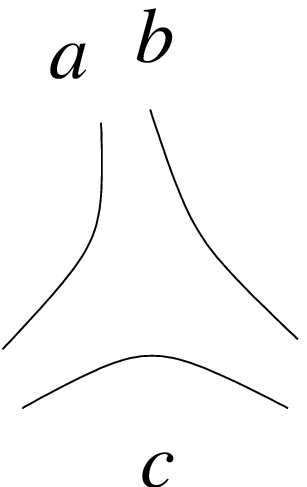}\right>
\qquad=\qquad \sum_{def} \frac{\dim_qd\,\dim_qe\,\dim_qf}{\theta_{abe}\theta_{acf}\theta_{bcd}\theta_{efd}}
\left<\epsfcenter{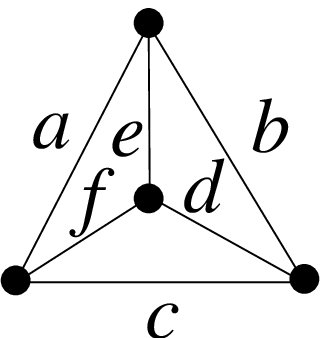}\right>\left<\epsfcenter{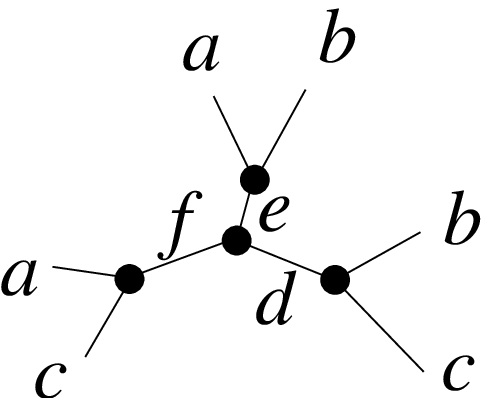}\right>
$$
$$
\left<\epsfcenter{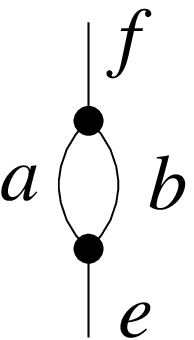}\right>\qquad=\qquad \delta_{fe}\frac{\theta_{abf}}{\dim_qf}\left<\quad\epsfcenter{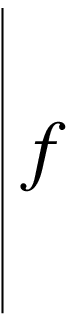}\right>
$$
In the case of a vertex which is four-valent or more, this can be reduced to the above cases by the identities such as $Z(M,\epsfcenter{fivevertex.eps})=Z(M,\epsfcenter{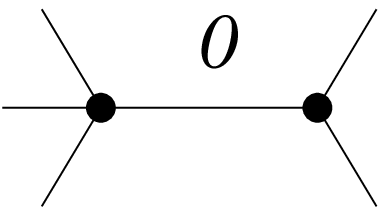})$, which follows from the fact that the chain-mail definitions of these two partition functions are the same. Finally, the overall normalisation factor is
$$ N^{-2-E-3C}N^{2C}=N^{-2-E-C}.$$
 
The result of applying the identites systematically in the case of a trivalent graph $\Gamma$ is the Kauffman bracket evaluation of the original labelled graph $\Gamma$, with a numerical weight factor which is equal to the Kauffman bracket evaluation of $\Gamma$ with all crossings reversed, multiplied by the inverse of a theta factor for each vertex.   
 This is best explained by giving the calculation for $\Gamma$ the graph with two edges given in the example above.
 \begin{multline*}
\sum _{jk} 
\frac 
{Z(S^3,\Gamma(j,k))} 
{Z(S^3)}
\frac{H_{jm}}{\dim_qj} 
\frac{H_{kn}}{\dim_qk}\\
= \frac{N^{-2-2-1}}{N^{-1}}\sum _{jk}\left<\, \epsfcenter{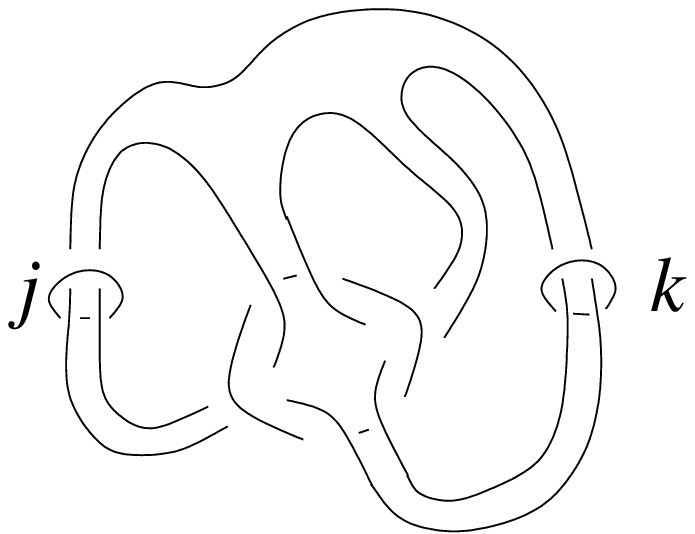}\right>H_{jm}H_{kn}\\
= N^{-1}\sum_{abcdl}
\frac{\dim_qa\dim_qb\dim_qc\dim_qd\dim_ql}{\theta_{abm}\theta_{adn}\theta_{bcn}\theta_{acl}\theta_{cdm}\theta_{lmn}^2}
\left<\epsfcenter{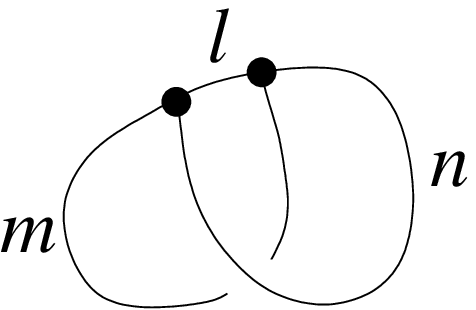}\right>\\
\left<\epsfcenter{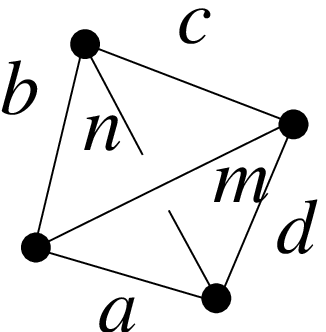}\right>
\left<\epsfcenter{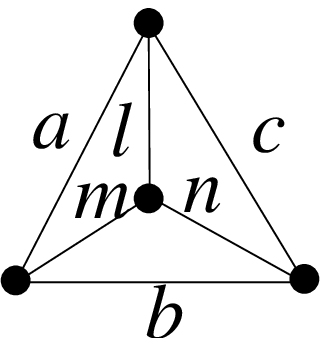}\right>
\left<\epsfcenter{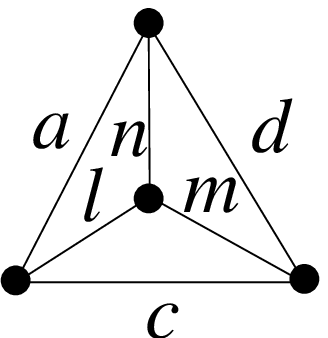}\right>\\
=\sum_l\frac{\dim_ql}{\theta_{lmn}^2}
\left<\epsfcenter{exampletriv.eps}\right>\left<\epsfcenter{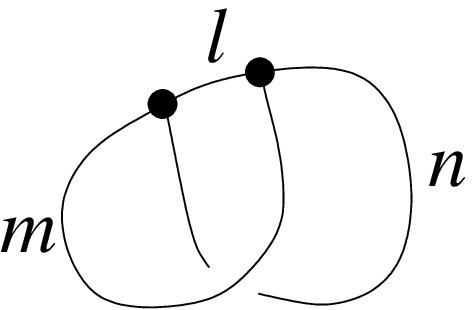}\right>
 = \left<\Gamma(m,n)\right>_R.
\end{multline*}  
This calculation uses the shadow world evaluation \cite{KR,KL}
\begin{multline*}\sum_{abcd}
\frac{\dim_qa\dim_qb\dim_qc\dim_qd}{\theta_{abm}\theta_{adn}\theta_{bcn}\theta_{acl}\theta_{cdm}}
\left<\epsfcenter{crosstet2.eps}\right>
\left<\epsfcenter{benztet2.eps}\right>
\left<\epsfcenter{benztet3.eps}\right>\\
=N  
 \left<\epsfcenter{exampletrivopp.eps}\right>.
 \end{multline*}
 
\section{Graphs in 3-manifolds}
 
In this section we give the definition of the graph invariants on manifolds followed by the proof of theorem \ref{theoremft}. The well-known Witten-Reshetikhin-Turaev invariant is used as the basis of the theory. In particular, many proofs are relatively simple because they exploit the known homeomorphism invariance of the Witten-Reshetikhin-Turaev invariant.

\subsection{The graph invariants on manifolds}
The Witten-Reshetikhin-Turaev invariant is an invariant of a labelled graph in an oriented closed 3-manifold $M$ which generalises the Kauffman bracket evaluation of a graph in $S^3$. This generalisation is required to define the $R$ invariant of theorem \ref{theoremft}, and also the 4-manifold invariants of later sections. 

Suppose $M$ is connected and $\gamma$ a framed graph in $M$ with each edge and each circle labelled with a spin, as in previous sections. Then $(M,\gamma)$ can be presented by the disjoint pair $\mu\cup\lambda\subset S^3$, where $\lambda$ is a framed labelled graph and $\mu$ a framed link, such that surgery on $\mu$ turns $(S^3,\lambda)$ into $(M,\gamma)$. Let $\kappa=\exp(i\pi(r-2)(3-2r)/4r)$, $m$ be the number of components of $\mu$, and $\sigma(\mu)$ the signature of the linking matrix of $\mu$. Then the invariant is defined by a sum over the possible spin labels for the components of the link $\mu$.

\begin{definition}[Witten-Reshetikhin-Turaev invariant]\label{wrtdef}
\begin{multline} \label{wrt}
Z_{WRT}(M,\gamma(j_1,j_2,\ldots,j_n))=\\
N^{-\frac{m+1}2}\kappa^{-\sigma(\mu)}
\sum_{i_1i_2\ldots i_m} \left<\mu\cup\lambda[i_1,i_2,\ldots, i_m,j_1,j_2,\ldots,j_n]\right>\dim_q i_1\dim_q i_2\ldots\dim_q i_m.
\end{multline}The definition is extended to manifolds with an arbitrary number of components by taking the product of (\ref{wrt}) for each component.
\end{definition}

Note that, in contrast to definition \ref{chainmail}, there is no factor of $\dim_q j$ for the labels which are not summed over. The invariant is sensitive to the orientation of $M$. Reversing the orientation of $M$ will be denoted $\overline M$.

Some elementary properties are
$$ Z_{WRT}(S^3)=N^{-1/2},\qquad Z_{WRT}(S^2\times S^1)=1,$$
$$ Z_{WRT}(S^3,\Gamma)=\left<\Gamma\right>Z_{WRT}(S^3),$$ 
$$ Z_{WRT}(\overline M,\Gamma)=\overline{Z_{WRT}(M,\Gamma)},$$
\begin{equation} \label{wrtproperties} Z_{WRT}(P\# Q)Z_{WRT}(S^3)= Z_{WRT}(P)Z_{WRT}(Q).\end{equation}
 The third equation follows from the fact that $N$ is real.

The relativistic invariant of a graph $\Gamma$ in $M$ is now defined in the same way as $<\ >_R$ was defined for $S^3$, in equation (\ref{reldef}). 

\begin{definition}[Relativistic invariant]
For a trivalent graph in an orientable 3-manifold $M$ 
\begin{equation}\label{reldefman}
Z_R(M,\Gamma) = 
\frac{\left| Z_{WRT}(M,\Gamma) \right| ^2} {\prod_{\text{vertices}}\Theta}.
\end{equation}
For the right-hand side of this definition it is necessary to pick a framing of the graph and an orientation of $M$, but the definition does not depend on these choices.
The definition is now extended to arbitrary graphs using the same relations as those for $<\ >_R$ following equation (\ref{reldef}). 
\end{definition}

The fact that this is a homeomorphism invariant follows from the fact that the 
Witten-Reshetikhin-Turaev invariant is, and the independence from framing follows by the same local calculations as in \cite{YO,Y}.

\subsection{Proof of Theorem \ref{theoremft}}\label{proofsection}
The case when the graph is empty reduces to
$$Z(M)=|Z_{WRT}(M)|^2$$
which was proved by Roberts \cite{RO}. We explain Roberts' method of proof and then extend it to deal with the case when the graph is not empty.

Roberts' method of proof was to show that in the case of a handle decomposition with one 0-handle and one 3-handle, performing surgery on the chain-mail diagram $\CH(D,\Phi)$ for oriented $M$ results in a manifold isomorphic to the connected sum of $M$ with its orientation reverse, $M\#\overline M$. Then he showed that $\sigma(\CH(D,\Phi))=0$, from which it follows that $Z(M)=N^{-1/2} Z_{WRT}(M\#\overline M)= |Z_{WRT}(M)|^2$.

We give this argument in more detail. Starting with a triangulation, we use lemma \ref{lemmaind} to convert the partition function to $Z_{CH}(M)$, for the corresponding handle decomposition with $h_k$ $k$-handles. In the chain-mail link $\CH(D,\Phi)$ the $\gamma_-$ circles have zero framing, because they bound embedded disks from $H_-$. 
It is easy to see that one can reduce the number of 3-handles to one without changing the value of the invariant by cancelling parallel $\gamma_-$ curves (this is part of theorem \ref{theoreminvt}, but one does not need to prove the entire theorem). 
The resulting link can be interpreted as a Kirby diagram for the 4-manifold $Y$ obtained from one 0-handle, 1-handles corresponding to the $\gamma_-$ curves and 2-handles corresponding to the $\gamma_+$ curves (and therefore the 1-handles of $M$).   According to the argument in \cite{GS}, sections 4.6.8 and 5.4, $Y=M_0\times [0,1]$, where $M_0$ is the 3-manifold $M$ with its 0-handles removed. 

To calculate the signature, we use the following principle \cite{RO}:

\begin{lemma} If $\mu$ is the link for a Kirby diagram for 4-manifold $W$, then $\sigma(\mu)=\sigma(W)$. \end{lemma}

If all the 1-handle curves are re-interpreted as 2-handle attaching curves, giving 4-manifold $V$, the result $\sigma(\mu)=\sigma(V)$ is well-known \cite{K}. However $W\cup_{\partial W=\partial V}\overline V$ is the boundary of a 5-manifold constructed from five-dimensional 2-handles; hence $\sigma(V)=\sigma(W)$.

Thus $\sigma(\CH(D,\Phi))=\sigma(Y)=0$. Also, standard arguments show that surgery on the link $\CH(D,\Phi)$ results in the 3-manifold $\partial Y$. But $\partial Y=M\#_{h_0}\overline M$, a generalised connected sum in which $h_0$ balls are removed before gluing the two summands together. This gives $Z(M)=N^{-h_0/2}Z_{WRT}(M\#_{h_0}\overline M)$.

Now starting with ${Z(M,\Gamma(j_1,j_2,\ldots,j_n))}$, the effect of the Fourier transform is to add a labelled framed link $L$ to the chain mail diagram for $M$. This is due to the identity
 $$ \left<\, \epsfcenter{line-j.eps}\, \right>H_{jk} 
=
\dim_q j \left<\,\epsfcenter{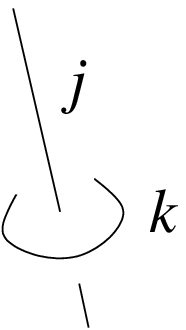} \right>$$
Therefore there is one additional circle labelled with $k_i$ corresponding to the $i$-th 1-handle in $\widetilde\Gamma$ in definition \ref{chainmail}. 

The Fourier transform is therefore
\begin{multline}\label{ftequation}
\sum _{{j_1j_2 \dots}  j_n} 
{Z(M,\Gamma(j_1,j_2,\ldots,j_n))} 
\frac{H_{j_1 k_1}}{\dim_qj_1} 
\frac{H_{j_2k_2}}{\dim_qj_2}\dots \frac{H_{j_n k_n}}{\dim_qj_n}\\
=N^{g-K-1} \sum _{j_{1}j_{2} \dots  j_K}
  \left<\left(\CH(D,\Phi)\cup L\right)[j_1,\ldots,j_K,k_1,\ldots k_n]\right>\;\prod_{i=1}^K \dim_q j_i \\
=
N^{-h_0/2} Z_{WRT}(M\#_{h_0}\overline M,L'(k_1,k_2,\ldots,k_n))
\end{multline}
where $L'$ is the link in $M$ obtained from the link $L$ in $S^3$ after surgery on the chain-mail diagram. The normalisation factors are $g-K-1=-h_0-h_2$, $m=h_1+h_2$ in equation (\ref{wrt}), and $h_3=1$; the Euler number $h_0-h_1+h_2-h_3=0$. 

The next step in the argument is to locate the components of the link $L'$. To do this, a more detailed description of the 3-manifold obtained by surgery on the original chain-mail link $\CH(D,\Phi)$ (before reducing the number of 3-handles to one) is given.  
Let $\Sigma=\partial H_-$ be the Heegaard surface. Surgery on the $\gamma_+$ curves in $\Sigma\times [0,1]$ gives $H_+\#_{h_0}{\overline H}_+$, and surgery on the  $\gamma_-$ curves in $\Sigma\times [0,1]$ gives ${\overline H}_-\#_{h_3}H_-$.

Consider the decomposition
$$S^3\cong H_-\cup\left(\Sigma\times [0,1]\right)\cup\left(\Sigma\times [0,1]\right)\cup(S^3\setminus \interior(H_-)),$$
using the convention that adjacent terms in the expression are glued by the obvious boundary components, and identifying $H_-$ with its image under $\Phi$.
The $\gamma_+$ curves lie in the first $ \Sigma\times [0,1]$ factor whilst the $\gamma_-$ curves lie in the second.
After surgery, $S^3$ becomes 
\begin{equation}\label{surgered}H_-\cup \left(H_+\#_{h_0} {\overline H}_+\right)\cup \left( {\overline H}_-\#_{h_3} H_-\right) \cup (S^3\setminus \interior(H_-))=M\#_{h_0}\overline M\#_{h_3} S^3.\end{equation} 

The surgery can be understood by the following argument due to Rourke \cite{RK}. To construct the result of surgery on a single framed curve $\gamma\subset \Sigma\times\{\frac12\}$, thicken $\gamma$ to an annulus $A\subset\Sigma\times\{\frac12\}$ (this is its framing), then glue a thickened disk $D^2\times [0,1]$ to $\Sigma\times[0,\frac12]$ along the attaching curve $\gamma$, using its thickening $A$. This results in a manifold $B_0$. Carry out a similar construction gluing a thickened disk to $\Sigma\times[\frac12,1]$ forming manifold $B_1$. Finally, $B_0$ and $B_1$ are glued along their two modified boundary components in the obvious way.

In fact each $\gamma_+$ curve gives rise to one 1-handle for $M$ and one 1-handle for $\overline M$ in the right-hand side of (\ref{surgered}), corresponding exactly to the original 1-handles of the manifold. The connected sum is removing a ball in $H_+$ for each of the original 0-handles. Exactly similar considerations apply to $H_-$ and its 2- and 3-handles.

  There is one circle $l$ of the link $L'$ for each of the 1-handles in the original decomposition of $M$. Before the surgery, the circle is an unknot linking the circle $\gamma_+$. Of course, it is not unique, since it is only determined up to ambient isotopy. Choose a simple (non-self-intersecting) arc $\lambda$ on the annulus $A$, defined above, which starts on one boundary circle of $A$ and ends on the other (see figure). 
$$\epsfbox{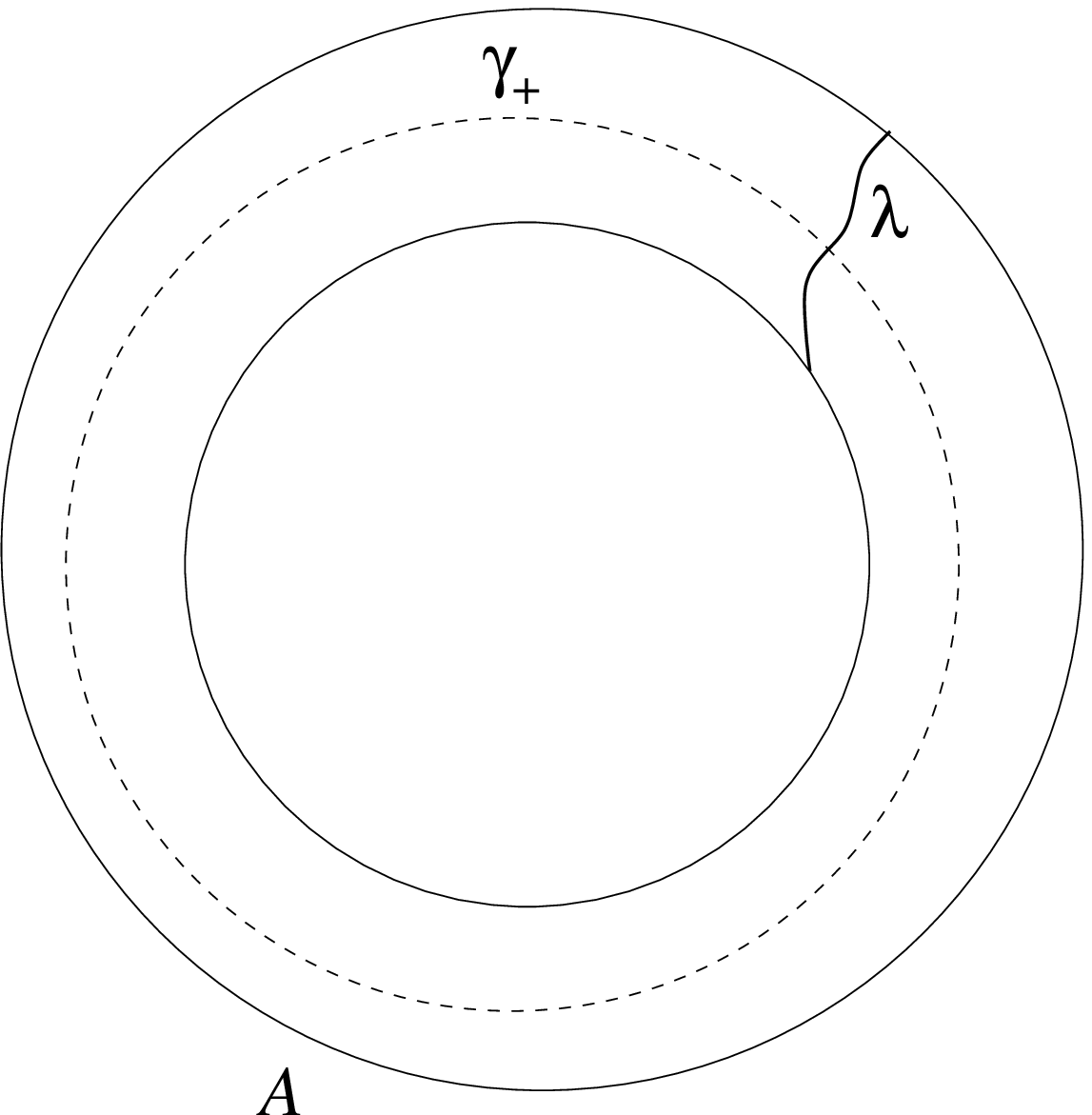}$$
This has a framing determined by thickening in $A$. One choice of the unknot $l$ is to take it to run along $\lambda$, deformed a small distance (in $\Sigma\times[0,1]$) to one side of $\gamma$, then back along $\lambda$, deformed to the other side of $\gamma$. After the surgery, the curve $l$ still runs along the two copies of $\lambda$, one in the factor $B_0$, the other in $B_1$ (it no longer need be deformed), as in the figure, which shows $H_+\#_{h_0}H_+$. Note that different choices of $\lambda$ (winding around $A$) correspond to different framings of the segment of $l$ in $M$, and the opposite framing of the segment of $l$ in $\overline M$. Hence the different choices of $\lambda$ all give the same framed curve $l$ up to isotopy.
$$ \epsfbox{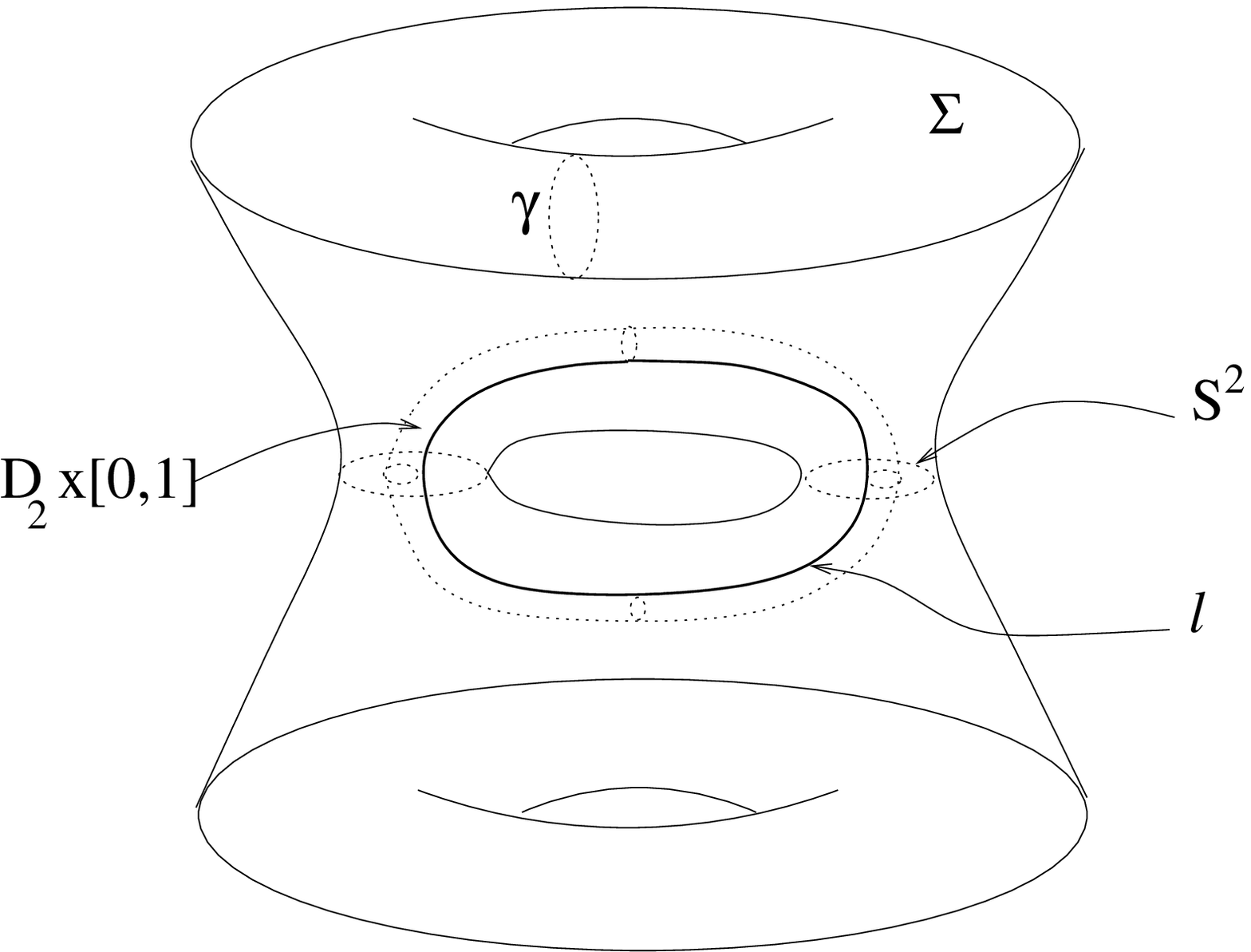}$$
It is clear that one can again reduce $h_3$ to one without changing this argument, thus giving the partition function in the last line of (\ref{ftequation}).

Finally, one can introduce an unlinked unknot in each $S^2$ connecting $M$ to $\overline M$, by multiplying with 
$$1=\frac1N\sum_j\dim_qj\left<\;\epsfcenter{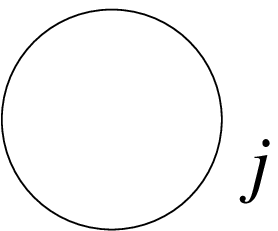}\;\right>.$$
  Due to the topology of $S^2$, this encircles the components of the link $L'$ which pass through it, and one can therefore apply the encircling lemma, shown here for the case of four strands
\begin{equation}\label{fourstrand}\sum_j\dim_qj \left<\epsfcenter{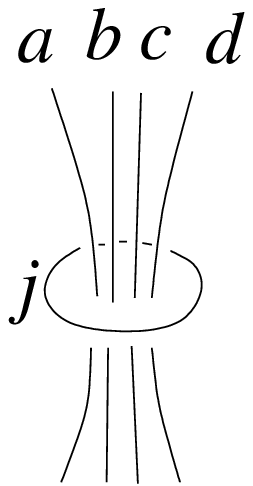}\right>
=N\sum_l
\theta_{abl}\theta_{cdl}
\;\dim_ql\;\left<\epsfcenter{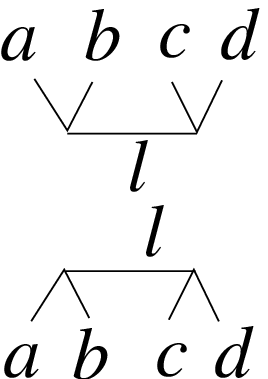}\right>.\end{equation}  
Then using the fact that $A\#_nB=A\#_{n-1}B\#(S^2\times S^1)$ and properties (\ref{wrtproperties}), one can remove the $\#_{h_0}$ in equation (\ref{ftequation}) at a cost of $N^{h_0/2}$. This gives 
$$N^{-h_0/2} Z_{WRT}(M\#_{h_0}\overline M,L'(k_1,k_2,\ldots,k_n))=
Z_R(M,\Gamma(k_1,k_2,\ldots,k_n)).$$
which proves the theorem.

Theorem \ref{theoreminvt} follows as an immediate corollary, since the relativistic invariant is known to be a homeomorphism invariant, using standard results about the Witten-Reshetikhin-Turaev invariant.

 \subsection{Alternative representation}
 One can also represent ${Z(M,\Gamma(j_1,j_2,\ldots,j_n))}$ as an invariant of the manifold  $M \#_{\widetilde\Gamma} \overline M$   obtained by removing the interior of the thickened graph $\widetilde\Gamma$ from both $M$ and $\overline M$ and gluing them together. This manifold has a labelled link $C(j_1,j_2,\ldots,j_n)$  given by the meridians of the thickened  edges of the graph $\Gamma$, with framings parallel to $\partial {\widetilde\Gamma}$, and the obvious colourings.  The relation is

\begin{theorem}\label{wrtobs}
$$Z(M,\Gamma(j_1,j_2,\ldots,j_n))=N^{-\frac{v+n}{2}}Z_{WRT}( M \#_{\widetilde\Gamma} \overline M, C(j_1,j_2,\ldots,j_n)   )  \prod_{i=1}^n \dim_q{j_i} ,$$
where $v$ is the number of vertices of $\Gamma$.
\end{theorem}

\noindent{\em Proof.} 
Let $D$ be a Heegaard diagram of a handle decomposition of $M$ with one 3-handle, as in the proof of Theorem \ref{theoremft}.
Write $\mathrm{CH }(D,\Phi)=\mu\cup\lambda$, with $\lambda$ the sublink corresponding to the edges of $\Gamma$, and $\mu$ its complement. From the same argument used in the proof of Theorem \ref{theoremft} it follows that the link $\mu$ is a surgery presentation of $M \#_{\n{\Gamma} \cup h_0} \overline{M}$ and the link $C$ is the link $\lambda$ after this surgery. If we see the 2-handle curves of  $\mu$   as dotted circles, then  $\mu$ is a Kirby diagram for  
$\left (M \setminus \interior(\n{\Gamma} \cup h_0) \right)\times I$. 
Thus $\sigma(\mu)=0$. 

Let $h_i$ be the number of $i$-handles of $M$. Then
\begin{align*}
 Z(M,\Gamma(j_1,j_2,\ldots,j_n)) &=Z_{\CH}(M,\widetilde\Gamma(j_1,j_2,\ldots,j_n))\\
 &=N^{g-K-1} \sum _{{j_{n+1}j_{n+2} \dots}  j_K}
  \left<\CH(D,\Phi)[j_1,\ldots,j_K]\right>\;\prod_{i=1}^K \dim_q j_i \\ 
&=N^{-(h_0+n)/2} Z_{WRT}(M\#_{h_0 \cup \n{\Gamma} }\overline M, C(j_1,j_2,\ldots,j_n)  ) \prod_{i=1}^n \dim_q{j_i}\\
 &=N^{-\frac{v+n}{2}}Z_{WRT}( M \#_{\widetilde\Gamma} \overline M, C(j_1,j_2,\ldots,j_n)   )  \prod_{i=1}^n \dim_q{j_i},
\end{align*}
by removing the 0-handles of $M$ which are not in $\widetilde\Gamma$.

\subsection{Definition for arbitrary graph}\label{3dInvariance}
It is now possible to determine the behaviour of the invariant $Z(M,\Gamma)$ under a subdivision of the graph $\Gamma$. Using the analogue of relation (\ref{twovalent}) for $Z_R$, one can use the inverse of the Fourier transform, together with the orthogonality
$$ \sum_k H_{jk}H_{kl}=N\delta_{jl}$$
to show that for any 2-valent vertex in a graph,
\begin{equation}\label{subdivide} Z(M, \quad\epsfpad{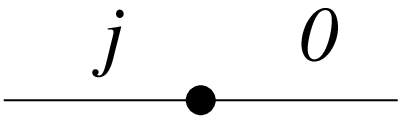}{0pt}{10pt}\quad)=\frac1N \;Z(M, \quad\epsfpad{line.eps}{0pt}{10pt}\quad),\end{equation}
providing $Z$ is defined for the graph obtained by removing the vertex.
Again, only the relevant portion of the graph is shown; it is assumed that the remainder of the graph is the same on both sides. 

In fact this result generalises: if the edges are labelled with $j$ and $k$
(rather than 0), then the partition function depends only on the product
$j\otimes k$ and dimension factors. More precisely:
\begin{multline}\label{subdivide2} \frac{1}{\dim_q j \dim_q k}Z(M, \quad\epsfpad{twovertex.eps}{0pt}{10pt}\quad)\\=\frac1N \sum_{a=|j-k|}^{j+k} \; \frac{1}{\dim_q a}Z(M, \quad\epsfpad{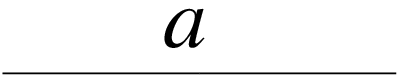}{0pt}{10pt}\quad),\end{multline}
 This identity can be  shown directly from the chain mail version  of the partition function, equation  (\ref{chainmail}), or Theorem \ref{wrtobs}. Alternatively, it is possible to use Theorem  \ref{theoremft} together with the inverse of the Fourier transform.

The definition of the invariant can now be extended to an arbitrary embedding of a graph in a closed 3-manifold. Let $G$ be an arbitrary graph as defined at the beginning of section \ref{rsn} (a compact 1-dimensional polyhedron), with a spin label for each edge or circle. Then there is a subdivision $\Gamma_1$ of $G$, and a triangulation of the 3-manifold $M$ such that $\Gamma_1$ is a simplicial complex and the embedding $\Gamma_1\to M$ is a simplicial map. Each edge or circle is subdivided into a number of edges of $\Gamma_1$. If $j$ is its spin label, then label one of the corresponding edges of $\Gamma_1$ with $j$ and the others with $0$. It is then possible to define the invariant of $G$ in terms of the invariant of $\Gamma_1$ given by the previous definition. Let $p$ be the number of additional vertices introduced by the subdivision of $G$. 

\begin{definition} \label{3dgeneral} $$Z(M,G(j_1,j_2,\ldots))=N^p Z(M,\Gamma_1(j_1,0,0,\ldots,j_2,0,\ldots)).$$
\end{definition}
This definition agrees with the previous definition in the case that $G$ is already a subset of the edges of some triangulation, due to theorem \ref{theoreminvt} and equation (\ref{subdivide}), and gives a homeomorphism invariant of the pair $(M,G)$. It also agrees with the chain-mail definition whenever $G$ is given by the cores of 0- and 1-handles for a handle decomposition. The results theorems \ref{theoremft} and \ref{theoremftcase}, and also equation (\ref{subdivide}) all hold with definition \ref{3dgeneral}.

\section{Four-manifolds}

Let $W$ be a closed connected oriented 4-manifold. The Crane-Yetter invariant of $W$ \cite{CKY} is defined by taking a triangulation of $W$ and defining a state-sum model on the triangulation. In this state sum, each triangle in the triangulation is labelled with a spin variable and the partition function $Z(W)\in \C$ is defined by summing the weights for each labelling over all values of the spin label for each triangle. The full definition is not given here, but an equivalent definition in terms of chain-mail is given below.

In a similar way to the previous case of 3-manifolds, one can consider defining an observable in the Crane-Yetter state sum by fixing the values of the spins on a subset of the triangles and summing over the remaining spin variables. Let $\Gamma$ be the 2-complex formed by this subset of triangles, and $\Gamma(j_1,j_2,\ldots,j_n)$ this 2-complex with the $i$-th triangle labelled by spin $j_i$. Then the state sum defines the observable
$$ Z\bigl(W,\Gamma(j_1,j_2,\ldots,j_n)\bigr)\in\C. $$

The principal tool will again be the chain-mail description of this observable in terms of a handle decomposition of $W$, generalising the description of the manifold invariant ($\Gamma=\emptyset$) by Roberts \cite{RO}. This definition is given first, followed by the equivalence with the Crane-Yetter definition for a triangulation.

A handle decomposition of $W$ with one 0-handle gives rise to a Kirby diagram, which is a link $L$ in $S^3$ with one component corresponding to each 1-handle or 2-handle of $W$ \cite{K}. The analogous observable for a handle decomposition is defined by the following formula. Let $\gamma$ be a sub-link of $L$ corresponding to a subset of the 2-handles, with the $i$-th component labelled by spin $j_i$. Order the components of $L$ so that the first $n$ components are the elements of $\gamma$. Let $h_i$ be the number of $i$-handles of $W$. Then the definition of the observable is
\begin{multline} \label{fourdch} Z_{CH}\bigl(W,\gamma(j_1,j_2,\ldots,j_n)\bigr)= \\
N^{-\frac12(h_1+h_2+h_3-h_4+1)}\sum_{j_{n+1}j_{n+2}\ldots j_K}\left<L[j_1,\ldots,j_n,\ldots,j_K]\right>\prod_{i=1}^K\dim_q j_i.
\end{multline}
 
The corresponding manifold invariant is 
$$Z(W)=\sum_{j_1\ldots j_n}Z_{CH}\bigl(W,\gamma(j_1,\ldots,j_n)\bigr)=\kappa^{\sigma(W)},$$
as can be seen by comparing the definition with \cite{RO}, using the fact that the number of components of $L$ is $h_1+h_2$ and the nullity of the linking matrix of $L$ is $h_3-h_4$. 

For each triangulation of $W$ we construct a handle decomposition which has one 0-handle.
The triangulation of $W$ gives rise to a canonical handle decomposition of $W$; however this needs to be modified. To give the correct correspondence with the Crane-Yetter state sum, the dual handle decomposition is used, so that the thickening of a vertex is a 4-handle, the thickening of an edge a 3-handle, etc. Then the number of 0-handles is reduced to one by replacing a simply connected tree containing all the 0-handles, and sufficient 1-handles to connect them, by a single 0-handle. This is now the correct handle decomposition of $W$.
 
 Then the relationship between these partition functions is:
\begin{lemma} \label{lemmafdch} Let $\gamma$ be the sublink of the Kirby diagram determined by the 2-handles corresponding to the triangles of $\Gamma$. Then
$$ Z\bigl(W,\Gamma(j_1,j_2,\ldots,j_n)\bigr)=Z_{CH}\bigl(W,\gamma(j_1,j_2,\ldots,j_n)\bigr).$$
\end{lemma}
This is proved simply by applying the relation (\ref{fourstrand}) to the component of $L$ corresponding to each tetrahedron in the triangulation in the definition of the right-hand side. This gives a sum of products of 15$j$-symbols and a normalisation factor of $N^{-\frac12\chi(W)}$ times the Crane-Yetter definition.

Now follows the proof of the main result given in the introduction which characterises the observable.

\subsection{Proof of Theorem \ref{theoremfd}}
The regular neighbourhood $\widetilde\Gamma$ of $\Gamma$ can be realised as the union of the 2- 3- and 4-handles which correspond to the 2- 1- and 0-simplexes of the 2-complex $\Gamma$. As above, these 2-handles can be taken to be the first $n$ in some ordering of the components of the link $L$ in the Kirby diagram.

According to lemma \ref{lemmafdch} and definition (\ref{fourdch}),
\begin{multline*}Z\bigl(W,\Gamma(j_1,j_2,\ldots,j_n)\bigr) = \\ 
 N^{-\frac12(h_1+h_2+h_3-h_4+1)}\sum_{j_{n+1}j_{n+2}\ldots j_K}\left<L[j_1,\ldots,j_n,\ldots,j_K]\right>\prod_{i=1}^K\dim_q j_i.
 \end{multline*}
  
The Kirby diagram omitting the first $n$ 2-handles describes the 4-manifold $W'$, the union of the 0-handle, all 1-handles, and the 2-handles not in $\widetilde\Gamma$. Thus
\begin{multline*}Z_{WRT}\bigl(\partial W',\gamma(j_1,j_2,\ldots,j_n)\bigr)=\\ N^{-\frac12(h_1+h_2-n+1)}\kappa^{-\sigma(W')} 
\sum_{j_{n+1}j_{n+2}\ldots j_K}\left<L[j_1,\ldots,j_n,\ldots,j_K]\right>\prod_{i=n+1}^K\dim_q j_i
\end{multline*}

Since $W'$ is obtained from $W\setminus\interior(\widetilde\Gamma)$ by successively removing 4-handles and then 3-handles, one can describe the effect of passing from  $\partial(W\setminus\interior(\widetilde\Gamma))$ to $\partial W'$ as follows.
 Firstly, removing each 4-handle adds a new $S^3$ component. Then, removing a 3-handle results in either joining two components together by the connected sum operation, or in forming the connected sum of one of the existing components with $S^1\times S^2$. 

Using the identities (\ref{wrtproperties}), this results in the relation
\begin{multline*}Z_{WRT}\bigl(\partial W',\gamma(j_1,j_2,\ldots,j_n)\bigr)=\\
N^{\frac12(h_3-h_4+\chi(\Gamma)-n)}
Z_{WRT}\bigl(\partial(W\setminus\interior(\widetilde\Gamma)),\gamma(j_1,j_2,\ldots,j_n)\bigr)
\end{multline*}
Finally, comparing these three equations, and using $\sigma(W')=\sigma(W\setminus\interior(\widetilde\Gamma))$, which
 follows from the fact that gluing 3- and 4-handles does not affect the signature of a 4-manifold \cite{W}, results in the statement of the theorem. 

\subsection{Examples}\label{sectionexamples}
According to theorem 2, the partition function for a 2-complex in a four-manifold $W$ depends only on a certain link in the boundary of a regular neighbourhood of the 2-complex. This involves a choice of a regular neighbourhood. All of these choices are related by a homeomorphism of $W$, and so the choice of regular neighbourhood does not affect the invariant. However it is not so easy to prove that the components of the link are well-behaved with respect to these homeomorphisms, and so a direct proof of the homeomorphism invariance of the partition function of a four-manifold with observables is missing. Nevertheless, in some specific cases it is easy to do the computations and see the invariance.

\subsubsection{Locally flat surfaces}\label{locallyflat}

Let $\S$ be a locally flat submanifold of an oriented  4-manifold $W$. Let $\n{\S}$ be a regular
neighbourhood of $\S$, thus $\n{\S}$ is a $D^2$-bundle over $\S$ with structure
group $S^1$.  See \cite{KOS} for  a discussion of these issues in the
smooth case. The piecewise linear version is treated in \cite{RS2}.

It follows that 
$\partial  \overline{\n{\S}}\cong  \partial (W \setminus 
\interior{\n{\S}}
)$  is an  $S^1$-bundle over
$\S$. If $\S$ is orientable, the  Euler class of this bundle is given by minus  the self intersection number of $\S$ with itself in
$W$.  This Euler class is called the normal Euler number of $\S$,
 and is always zero  for   oriented  surfaces in $S^4$, since
 $H_2(S^4)=0$. However, the same is not  true for  non-orientable surfaces \cite{CKS}.

Explicit examples where $\S$ is a locally flat surface can be computed easily. 
 In the case that $\S$ is an $S^2$, the observable is
 \begin{multline*}Z_{WRT}\bigl(\partial\overline{\widetilde\S},\gamma(j_1,\ldots j_n)\bigr)
 \kappa^{(\sigma(W)-\sigma(\widetilde\S))}
 N^{1-n}\dim_qj_1\ldots\dim_q j_n
 \\=
 \kappa^{\sigma(W)}N^{-n}\dim_qj_1\ldots\dim_q j_n \sum_i\dim_qi
 \quad \left<\;\epsfcenter{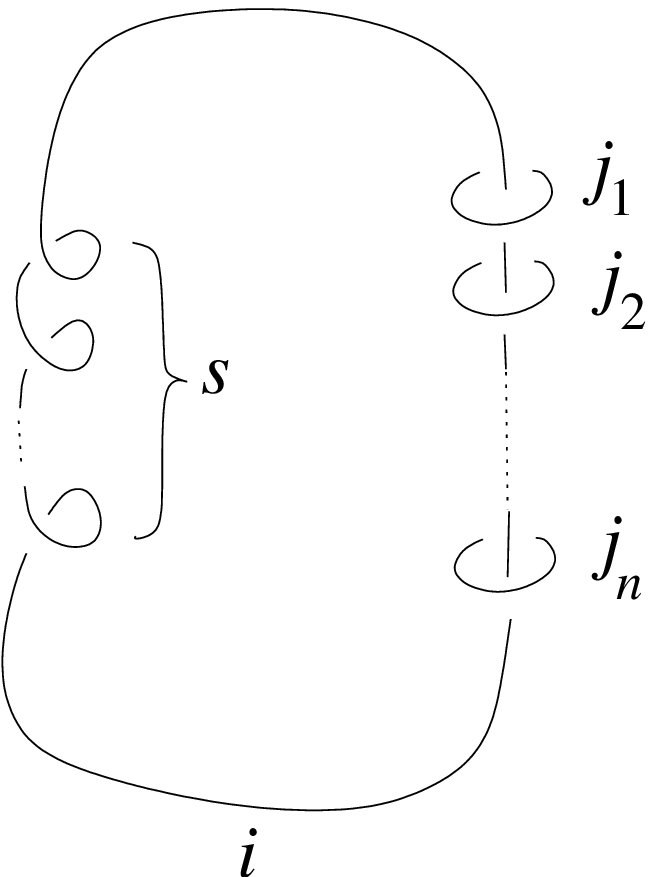}\;\right>,\end{multline*} 
 where $s$ is the normal Euler number of $\S$.

 In the case that $\S$ is an $S^1\times S^1$, the formula is 
 $$\kappa^{\sigma(W)}N^{-n-2}\dim_qj_1\ldots\dim_q j_n \sum_{abc}\dim_qa\dim_qb\dim_qc
 \quad \left<\;\epsfcenter{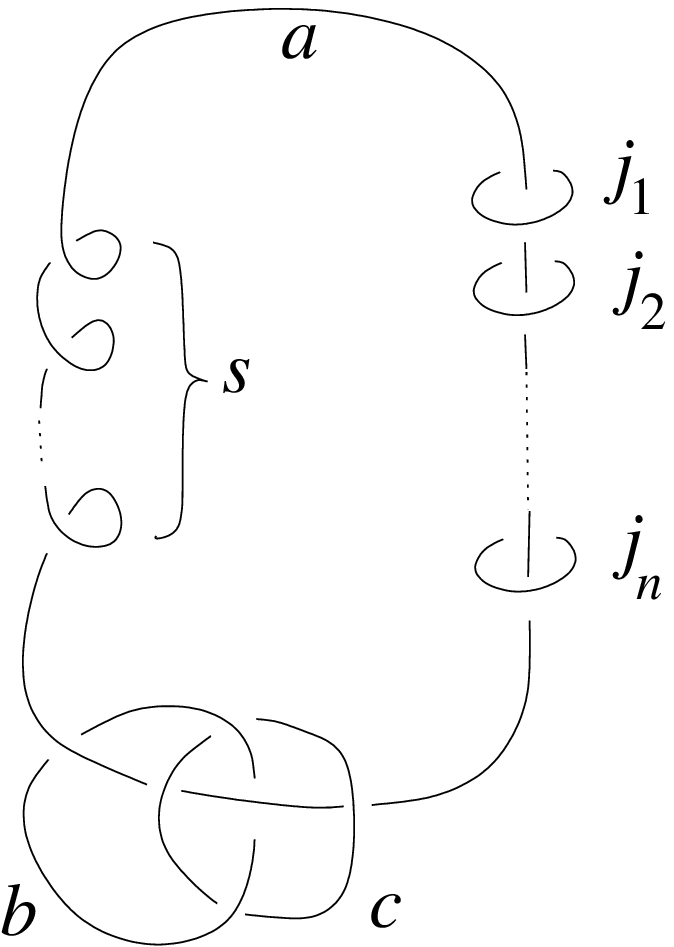}\;\right>,
 $$
where as before $s$ is the normal Euler number of $\S$. 
 The curves labelled with $b$ and $c$ are the 1-handles whilst the curve labelled by $a$ is the attaching curve for the 2-handle of $\widetilde \S$.
 
From these examples, it can be seen that $Z(W,\Sigma(j_1,\ldots,j_n))$ depends only on the tensor product of the representations colouring $\Sigma$, and it can easily be normalised to give a triangulation-independent quantity. However, the homeomorphism invariant obtained depends only on the Euler number of the normal bundle of the surface.

\subsubsection{Trivial labelling}

 One can use the statement of theorem \ref{theoremfd} to extend the definition of $Z$ to include the case of oriented 4-manifolds with boundary. For example, without observables this gives
 $$ Z(W)=Z_{WRT}(\partial W)\kappa^{\sigma(W)}$$
 This has the property that if the closed three-manifold $M=W_1\cap W_2$ is a union of components of the boundary of both $W_1$ and $W_2$, then 
 $$ Z(W_1)Z(W_2)=Z(M) Z(W_1\cup W_2),$$
with $Z(M)$ being the Turaev-Viro invariant of $M$.
Also, for any closed 3-manifold $M$,
$$ Z(M\times I)=Z(M).$$ 
This follows from the fact that $\sigma(M\times I)=0$.
 
These invariants are in fact obtained from the case of the observables on a 2-complex in a manifold without boundary, by taking all of the labels to be $0$. This has the effect of removing a regular neighbourhood of the 2-complex.

\begin{theorem} \label{twocx}
Let $\Gamma$ be a 2-dimensional polyhedron. Consider a piecewise linear embedding of
$\Gamma$ into an oriented 4-dimensional manifold $W$.  Choose  a triangulation of
$\Gamma$ so that this triangulation can be extended to a triangulation of
$W$ (this is always possible). Let $n$ be the number of 2-simplexes of 
this 
triangulation of $\Gamma$. The quantity $N^nZ(W,\Gamma(0,0,\ldots,0))$
is a homeomorphism invariant of $\Gamma$ in $W$. In fact
\begin{align*}
N^nZ(W,\Gamma(0,0,\ldots,0))&= Z_{WRT}\bigl(\partial(W\setminus
\interior(\widetilde\Gamma))\bigr)
\kappa^{\sigma(W\setminus\interior(\widetilde\Gamma))} N^{\frac{\chi(\Gamma)}{2}}\\
                 &= Z(W\setminus\interior(\widetilde\Gamma)) N^{\frac{\chi(\Gamma)}{2}}.
\end{align*}
\end{theorem}

As an example, consider a graph $\gamma\subset S^3$. The suspension $S\gamma\subset S^4$ is defined by taking the join with $S^0=\{a,b\}$. If $S^3$ is triangulated and $\gamma$ forms a subset of the edges, then the suspension gives a triangulation of $S^4$ and the 2-complex $S\gamma$ which determined by a subset of the two-simplexes. Thus there is a Crane-Yetter partition function. With trivial labellings for the two-complex, theorem \ref{twocx} gives the homeomorphism invariant of $S\gamma$ in $S^4$
 \begin{equation*}N^n Z\bigl(S^4, S\gamma(0,0,\ldots,0)\bigr)
\end{equation*}
where $n$ is the number of 
triangles of $S\gamma$. This invariant can be calculated from the three-dimensional invariant by the identity
 \begin{equation}N^n Z\bigl(S^4, S\gamma(0,0,\ldots,0)\bigr)=
N^{1+m}Z\bigl(S^3,\gamma(0,0,\ldots,0)\bigr),
\end{equation}
where $m$ is the number of edges of $\gamma$.
This identity is proved as follows. Let $U_1$ be a regular neighbourhood of $a$ in $S^4$ and $U_2$ a regular neighbourhood of $b$. Then $S^4$ is homeomorphic to $U_1\cup (S^3\times I)\cup U_2$, and a regular neighbourhood of $S\gamma$ is given by $U_1\cup (\widetilde\gamma\times I)\cup U_2$, 
with $\widetilde{\gamma}$ a regular neighbourhood of $\gamma$ in $S^3$; cf. the proof of theorem 2(a) of Cencelj, Repov\v{s} and Skopenkov \cite{CRS} for a similar construction.

Hence $S^4\setminus\interior( \widetilde {S\gamma})$ is homeomorphic to $(S^3\setminus\interior(\widetilde\gamma))\times I$, and has boundary $S^3\#_{\widetilde\gamma}\overline {S^3}$, as in the proof of theorem \ref{theoremft}. Thus
\begin{align*}N^n Z\bigl(S^4, S\gamma(0,0,\ldots,0)\bigr) &=
Z(S^4\setminus\interior(\widetilde{S\gamma}))N^{\chi(S\gamma)/2}\\
&=Z\bigl( (S^3\setminus \interior(\widetilde\gamma))\times I\bigr)N^{(2-v+m)/2}\\
&=Z_{WRT}\bigl(S^3\#_{\widetilde\gamma}\overline {S^3}\bigr)N^{(2-v+m)/2}\\
&=N^{1+m}Z\bigl(S^3,\gamma(0,0,\ldots,0)\bigr),
\end{align*}
the last equality using theorem \ref{wrtobs}.
This example shows that the observables are non-trivial and can be used to detect some non locally-flat embeddings of $S^2$ in $S^4$, by taking $\gamma$ to be a knot. Then $S\gamma$ is non locally-flat if the knot is non-trivial\cite{CRS}, and the invariant is non-trivial whenever $Z_{WRT}(S^3 \#_{\widetilde\gamma} \overline{S^3}) \neq 1$.

\subsection*{Acknowledgement}
Thanks are due to a conversation with Du\v{s}an Repov\v{s}. JFM was financed by  Funda\c{c}\~{a}o para a Ci\^{e}ncia e
Tecnologia (Portugal), post-doctoral grant number SFRH/BPD/17552/2004, and part of the research project POCTI/MAT/60352/2004 (``Quantum Topology''), also financed by the FCT.

\end{document}